\documentclass[mathematics,article,submit,moreauthors,pdftex]{mdpi}
\firstpage{1} 
\pubvolume{1}
\issuenum{1}
\articlenumber{0}
\pubyear{2021}
\copyrightyear{2020}
\datereceived{} 
\dateaccepted{} 
\datepublished{} 
\hreflink{https://doi.org/} 

\usepackage{mathtools}
\usepackage{bbm}
\usepackage{amsfonts}
\usepackage{amssymb}

\def\Bila{\mathbf{B}}
\def\S{\mathbf{S}}

\def\P{\mathcal{P}}
\def\J{\mathcal{J}}
\def\K{\mathcal{K}}
\def\D{\mathcal{D}}
\def\R{\mathbb{R}}
\def\DD{\mathbb{D}}

\def\A{\mathcal{A}}
\def\I{\mathcal{I}}
\def\F{\mathcal{F}}
\def\bu{\textbf{u}}
\newcommand {\nor} [1]{\parallel #1 \parallel}

\newcommand{\bequ}{\begin{equation}}
\newcommand{\eequ}{\end{equation}}
\newcommand{\bequd}{\begin{eqnarray*}}
\newcommand{\eequd}{\end{eqnarray*}}

\definecolor{Darkgreen}{rgb}{0,0.5,0}
\definecolor{Purple}{rgb}{0.5,0,0.5}
\definecolor{Orange}{rgb}{1,0.65,0}

\newcommand{\sz}[1]{\textcolor{black}{#1}}
\newcommand{\re}[1]{\textcolor{black}{#1}}
\newcommand{\ds}[1]{\textcolor{black}{#1}}
\newcommand{\khi}[1]{\textcolor{black}{#1}}

\Title{Mathematical modelling of glioblastomas invasion within the brain: a 3D multi-scale moving-boundary approach}

\TitleCitation{Mathematical modelling of glioblastomas within the brain: a 3D multi-scale moving-boundary approach}

\Author{Szabolcs Suveges $^{1,a}$, Kismet Hossain-Ibrahim$^{2,3}$, J. Douglas Steele$^{4}$, Raluca Eftimie $^{5}$ and Dumitru Trucu $^{1,b,}$*}

\AuthorNames{Szabolcs Suveges, Kismet Hossain-Ibrahim, J. Douglas Steele, Raluca Eftimie and Dumitru Trucu}

\AuthorCitation{Suveges, S.; Hossain-Ibrahim, K.; Steele, J.D.; Eftimie, R.; Trucu, D.}

\address{%
$^{1}$ \quad Division of Mathematics, University of Dundee, Dundee DD1 4HN, UK; $^{\textit{\textbf{a}}}$ssuveges@dundee.ac.uk; $^{\textit{\textbf{b}}}$trucu@maths.dundee.ac.uk\\
$^{2}$ \quad Division of Cellular and Molecular Medicine, School of Medicine, University of Dundee, Dundee, UK; \khi{kismet.ibrahim@nhs.scot}\\
$^{3}$ \quad Department of Neurosurgery, Ninewells Hospital and Medical School, NHS Tayside, Dundee, UK; kismet.ibrahim@nhs.scot\\
$^{4}$ \quad Division of Imaging Science and Technology, Medical School, University of Dundee, Dundee, UK; d.steele@dundee.ac.uk\\
$^{5}$ \quad Laboratoire Math\'{e}matiques de Besan\c{c}on, UMR - CNRS 6623, Universit\'{e} de Bourgogne Franche-Comt\'{e}, 16 Route de Gray, Besan\c{c}on, France; raluca.eftimie@univ-fcomte.fr}

\corres{Correspondence: trucu@maths.dundee.ac.uk}

\abstract{
\re{Brain-related experiments are limited by nature, and so biological insights are often restricted or absent. This is particularly problematic in the context of brain cancers, which have very poor survival rates. To generate and test new biological hypotheses, researchers started using mathematical models that can simulate tumour evolution. However, most of these models focus on single-scale 2D cell dynamics, and cannot capture the complex multi-scale tumour invasion patterns in 3D brains. A particular role in these invasion patterns is likely played by the distribution of micro-fibres. To investigate explicitly the role of brain micro-fibres in the 3D invading tumours, in this study we extend a previously-introduced 2D multi-scale moving-boundary framework to take into account 3D multi-scale tumour dynamics. \ds{T1 weighted and DTI scans are} used as initial conditions for our model, and to parametrise the diffusion tensor. Numerical results show that including an anisotropic diffusion term may lead in some cases (for specific micro-fibre distributions) to significant changes in tumour morphology, while in other cases it has no effect. This may be caused by the underlying brain structure and its microscopic fibre representation, which seems to influence cancer-invasion patterns through the underlying cell-adhesion process that overshadows the diffusion process.}
}


\begin{document}
\section{Introduction}
Glioblastoma multiforme is a highly invasive and aggressive type of brain tumour, typically with poor \re{patient} prognosis \cite{Burri2018,Davis2016,Klopfenstein2019,Louis2007,Meneceur2020,Preusser2011,Sottoriva2013} (\khi{median} survival rate is less than 1 year \cite{Brodbelt2015}). \re{These tumours} arise from abnormal glial cells located in the central nervous system, and shortly after their appearance they invade the surrounding tissues in a heterogeneous fashion. This \re{heterogeneous invasion pattern} leads to tumours whose outer edges are difficult or impossible to determine with current imaging technologies, including for instance \emph{magnetic resonance imaging} (MRI) and \emph{diffusion tensor imaging} (DTI), both of which \khi{measure} the diffusion of water molecules and enable the study of brain structures.

Due to the limited experimental \re{approaches} that one can \re{use to study} the brain, \re{researchers have started using} mathematical models \re{to} provide certain \re{biological} insights that otherwise would be difficult to obtain \re{experimentally}. \re{Such models can help predict how tumours grow for specific patients, aiding clinicians in decision-making, or they can help test and provide new hypotheses about potential anti-tumour treatments. 
}
Mathematical modelling of tumours has seen significant advances over the last few decades, which broadened our understanding of tumour dynamics and how cells interact with their environment \cite{Anderson_et_al_2000,Anderson2009,Anderson_2005,Basanta2008,Basanta2011,Boettger2012,Chaplain_Lolas_2005,Chaplain2006,Deakin_Chaplain_2013,Deisboeck_et_al_2011,Domschke_et_al_2014,Dumitru_et_al_2013,Hatzikirou2010,Kiran_2009,KnutsdottirPalssonKeshet2014,Macklin_et_al_2009,MahlbacherLowengrubFriesboes2018,Shuttleworth_2019,Shuttleworth2020a,Shuttleworth2020b,Suveges_2020,Suveges_2021,szymanska_08,Tektonidis2011,Xu_2016}. Although \re{the majority of these} models do not restrict themselves to a specific tumour type and rather focus on general tumours, there are some \re{models} that focus on the evolution of gliomas within the brain \cite{Alfonso2016,Engwer2014,Hunt2016,Painter2013,Scribner2014,Swanson2000,Swanson2007,Swanson2011,Sykova2008}. Recently, models also started to incorporate the structure of the brain, by including MRI and DTI \ds{scans} \cite{Clatz2005,Cobzas2009,Engwer2014,Hunt2016,Jbabdi2005,Konukoglu2010,Painter2013}. Even though these images are generated in 3D, most of these models are simulating the tumour growth in 2D and only a few of them are 3D models \cite{Clatz2005,Suarez2012,Yan2017}. Moreover, the majority of \re{published models focus} on tumour progression only on one spatio-temporal scale. However, \re{tumour progression is characterised} by various biological processes occurring on different scales, and \re{thus} their effects on the overall tumour dynamics cannot be neglected. Hence, recent efforts have been made to establish new multi-scale frameworks for tumour progression \cite{Engwer2014,Hunt2016,Painter2013,Peng2016,Shuttleworth_2019,Shuttleworth2020a,Shuttleworth2020b,Suveges_2020,Suveges_2021}, which were able to capture some of these multi-scale underlying biological processes usually involving the extracellular matrix (ECM). 

In this paper, we extend the general \re{2D} multi-scale moving-boundary modelling framework introduced in \cite{Dumitru_et_al_2013,Shuttleworth_2019} to capture the invasion of glioblastomas within a 3D fibrous brain environment. To this end, we incorporate the information provided by both the \ds{T1 weighted and DTI scans} into our \re{multi-scale} framework and use the resulting model to simulate \re{numerically} the growth of 3D gliomas within the brain. \re{We focus on a few cases showing tumour growth in different regions in the brain, with different distributions of grey/white matter densities, which leads to different tumour invasion patterns.}

The paper is organised as follows. First, we formulate our extended multi-scale moving boundary framework in Section~\ref{sec:Modelling}. Then, following a brief description of the numerical methods, we present the computational simulation results in Section~\ref{sec:Results}. Finally, in Section~\ref{sec:Discussion} we summarise and discuss these results.

\section{Multi-Scale Modelling of the Tumour Dynamics}\label{sec:Modelling}
To model the evolution of glioblastomas within a 3-dimensional brain, we employ a multi-scale moving boundary model that was initially introduced in \cite{Dumitru_et_al_2013} and later expanded in several other works \cite{Peng2016,Shuttleworth_2019,Shuttleworth2020a,Shuttleworth2020b,Suveges_2020,Suveges_2021}. To account for the brain's structure, we aim to use 3D \ds{T1 weighted and DTI scans} that ultimately influence the migration of the cancer cells as well as affect both micro-scale dynamics. Hence, here we aim to explore the impact of the brain structure on the interlinked macro-scale and micro-scale tumour dynamics.

\subsection{Macro-Scale Dynamics}
\begin{figure}[t!]
\centering
\includegraphics[width=10.5 cm]{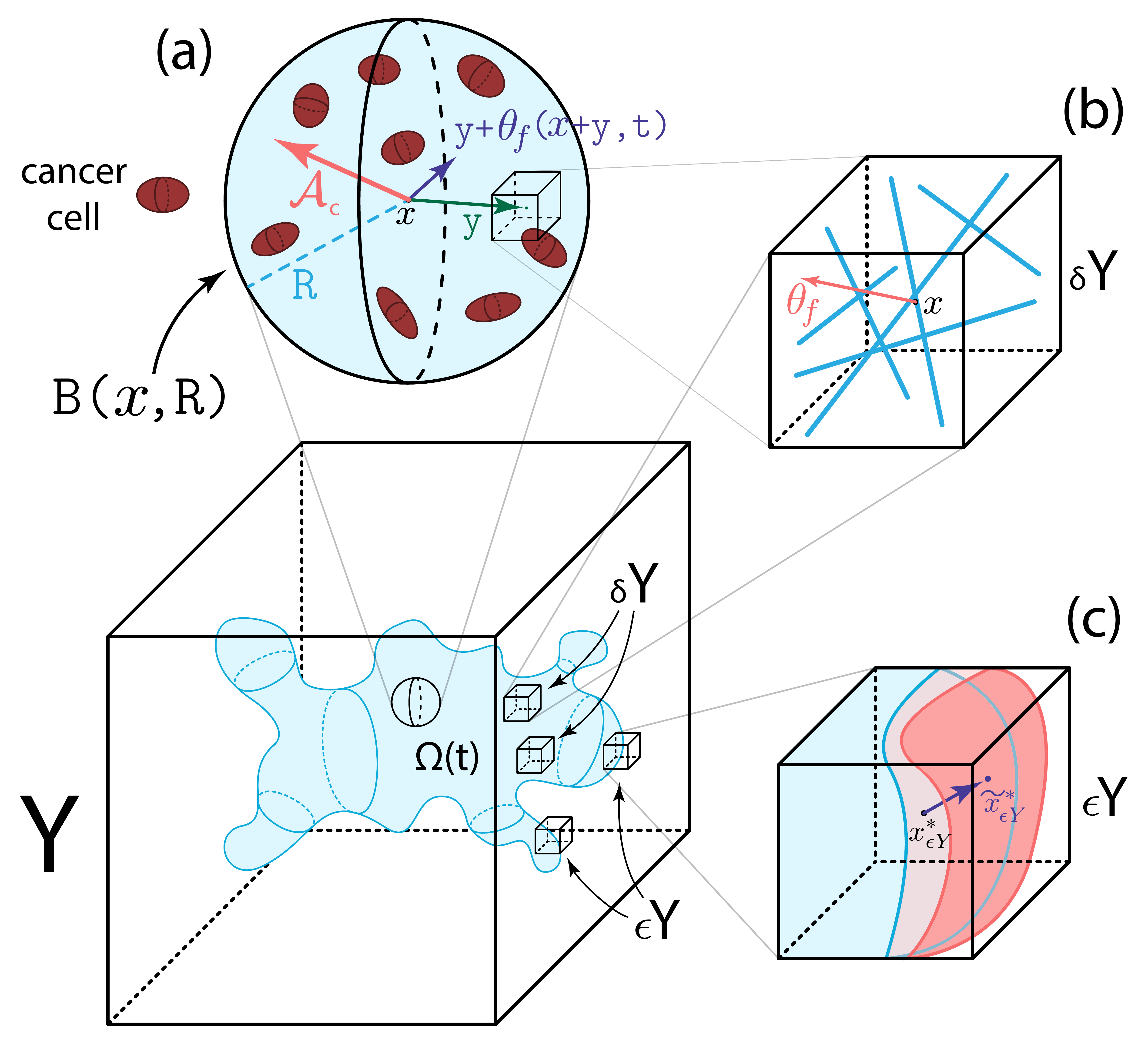}
\caption{\sz{Schematics of the multi-scale model. (a) Illustration of the sensing region $\Bila(x, R)$, the two vectors $y$ and $y + \theta_{f}(y+x,t)$ and the overall travelling direction $\A_{c}$. (b) An example of a fibre micro-domain $\delta Y(x)$ containing fibres (blue lines) that induces an overall orientation $\theta_{f}(x,t)$ for $\delta Y(x)$. (c) An example of a boundary micro-domain $\epsilon Y(x)$ where the blue volume represents the tumour volume at the current time-step with boundary point $x$ and the red volume represents the evolved tumour at next time-step with shifted boundary point $x_{\epsilon Y}^{*}$.}\label{fig:Model_Schematics}}
\end{figure}

Since in this work we extend the \re{2-dimensional (2D)} modelling framework introduced in \cite{Dumitru_et_al_2013,Shuttleworth_2019}, we begin by \re{describing briefly} some of the key features \re{of this framework} and \re{by giving} a few useful notations. First, we denote \re{by $\Omega(t)$} the expanding 3-dimensional \re{(3D)} tumour region that progresses over the time interval $[0,T]$ within a maximal tissue cube $Y \subset \R^{N}$ with \re{$N=3$}, \emph{i.e.,} $\Omega(t) \subset Y,\;\forall t \in [0,T]$; \re{see also Figure \ref{fig:Model_Schematics}.}
Then at any macro-scale spatio-temporal point $(x,t) \in Y \times [0,T]$ we consider a cancer cell population $c(x,t)$ embedded within a two-phase ECM, consisting of the non-fibre $l(x,t)$ and fibre $F(x,t)$ ECM phases \cite{Shuttleworth_2019,Shuttleworth2020a,Shuttleworth2020b,Suveges_2020,Suveges_2021}. On the one hand, the fibre ECM phase accounts for all major fibrous proteins such as collagen and fibronectin, whose micro-scale distribution induces the spatial orientation of the ECM fibres. Hence, the macro-scale spatio-temporal distribution of the ECM fibres \re{is} represented by an oriented vector field $\theta_{f}(x,t)$ that describes their spatial bias, as well as by $F(x,t):=\nor{\theta_{f}(x,t)}$ which denote the amount of fibres at a macro-scale point $(x,t)$ \cite{Shuttleworth_2019,Shuttleworth2020a,Shuttleworth2020b,Suveges_2020,Suveges_2021}. On the other hand, in the non-fibre ECM phase we bundle together every other ECM constituent such as non-fibrous proteins (for example amyloid fibrils), enzymes, polysaccharides and extracellular $Ca^{2+}$ ions \cite{Shuttleworth_2019,Shuttleworth2020a,Shuttleworth2020b,Suveges_2020,Suveges_2021}. 
\re{Furthermore, in this new modelling study we incorporate the structure of the brain by extracting data from \ds{DTI and T1 weighted brain scans}, and then using this data to parametrise the model. Specifically, }we denote by $\DD_{Water}(x)$ the water diffusion tensor that is induced by the DTI \ds{scan}. Also, we denote by $w(x)$ the white matter density and by $g(x)$ the grey matter density, both of which are extracted from the \ds{T1 weighted image}. Finally, for \re{compact writing}, we denote \re{by $\bu$} the global tumour vector described as
\begin{linenomath}\begin{equation*}
	\bu := (c(x,t), l(x,t), F(x,t))^{\intercal},
\end{equation*}\end{linenomath}
\re{and} by $\rho(\bu)$ the total space occupied that is defined as
\begin{linenomath}\begin{equation*}
	\rho(\bu) := c(x,t) + l(x,t) + F(x,t).
\end{equation*}\end{linenomath}

\subsubsection{Cancer cell population dynamics}
\re{To describe the spatio-temporal evolution of} the cancer cell population $c(x,t)$, \re{we first assume a logistic-type growth with rate $\mu$ \cite{Laird1964,Laird1965,Tjorve2017,Shuttleworth_2019,Shuttleworth2020a,Shuttleworth2020b,Suveges_2020,Suveges_2021}. For the movement of this cell population, we use} the structure of the brain by taking into account both the \ds{T1 weighted and DTI scans} \ds{(from the IXI Dataset \cite{IXIDataset})}, as well as the various adhesion mediated processes \cite{Chen_2011,Condeelis2006,Huda2018,Petrie2009,Weiger2013,Wu3949}. Hence, the spatio-temporal dynamics of the cancer cell population is described by
\begin{linenomath}\begin{equation}
	\dfrac{\partial c}{\partial t} = \underbrace{\nabla \nabla :\big[ \DD_{T}(x)c \big]}_{\substack{\text{fully anisotropic} \\ \text{diffusion}}} - \underbrace{\nabla \cdot \big[ c \A_{c}(x,t,\bu,\theta_{f}) \big]}_{\text{adhesion processes}} + \underbrace{\mu c \big[ 1 - \rho(\bu) \big]^{+}}_{\substack{\text{logistic-type} \\ \text{proliferation}}}.
	\label{Cancer_Cell_Dynamics}
\end{equation}\end{linenomath}
Here, the operator $\nabla \nabla :$ denotes the full second order derivative \cite{Engwer2014}, \emph{i.e.,} it is defined as
\begin{linenomath}\begin{equation*}
	\nabla \nabla :\big[ \DD_{T}(x)c \big] := \sum_{i,j=1}^{N} \dfrac{\partial}{\partial x_{i}} \dfrac{\partial}{\partial x_{j}} \Big( \DD_{i,j} c \Big), \qquad N=3,
\end{equation*}\end{linenomath}
with $\DD_{i,j}$ denoting the components of the tumour diffusion tensor $\DD_{T}$.
Since classical diffusion models with constant coefficient cannot capture any directional cues, \re{as those} provided by the DTI data, in Eq. \eqref{Cancer_Cell_Dynamics} we use a tensor model (involving a fully anisotropic diffusion term) that is able to incorporate the anisotropic nature of the cancer cell movement. These tensor models were proposed \re{in} \cite{Basser1992,Basser1993,Basser1994a,Basser1994b} and \re{have} been used to mathematically model the gliomas within the brain\re{;} see for instance \cite{Engwer2014,Hunt2016,Painter2013}. The \re{main} idea \re{of this approach} is to use the measured water diffusivity in the structured, fibrous brain environment characterised by a symmetric water diffusion tensor
\begin{linenomath}\begin{equation}
	\DD_{Water}(x) =
	\begin{bmatrix}
    	d_{xx}(x) & d_{xy}(x) & d_{xz}(x) \\
    	d_{xy}(x) & d_{yy}(x) & d_{yz}(x) \\
    	d_{xz}(x) & d_{yz}(x) & d_{zz}(x)
    \end{bmatrix},
    \label{Water_Diffusion_Tensor}
\end{equation}\end{linenomath}
and appropriately construct a macroscopic diffusion tensor for the cancer cell population. Since this water tensor \eqref{Water_Diffusion_Tensor} is assumed to be symmetric, it can be diagonalised. Denoting its eigenvalues by $\lambda_{1}(x)\geq \cdots \geq \lambda_{N}(x)$ and the associated eigenvectors by $\phi_{1}(x), \cdots \phi_{N}(x)$, we follow \cite{Hillen2017,Mardia2000,Painter2013} and construct the \re{3D} tumour diffusion tensor as
\begin{linenomath}\begin{equation}
\begin{split}
	\DD_{T}(x) := D_{c} \; D_{WG}(x) \bigg[ & \bigg( r + (1-r)\bigg( \dfrac{\coth k(x)}{k(x)} - \dfrac{1}{k(x)^2} \bigg) \bigg) \re{I_{3}} \\
	& + (1-r)\bigg( 1 - \dfrac{3\coth k(x)}{k(x)} + \dfrac{3}{k(x)^2} \bigg) \phi_{1}(x) \phi_{1}^{\intercal}(x) \bigg].
	\label{Tumour_Diffusion_Tensor}
\end{split}
\end{equation}\end{linenomath}
\re{Here} $D_{c} > 0$ is the diffusion coefficient, $r \in [0, 1]$ specifies the degree of isotropic diffusion, $\re{I_{3}}$ denotes the $3 \times 3$ identity matrix, $k(x)$ is given by
\begin{linenomath}\begin{equation*}
	k(x) := \K_{FA} FA(x),
\end{equation*}\end{linenomath}
with $\K_{FA} \geq 0$ being a proportionality constant measuring the sensitivity of the cells to the environments' direction, and $FA(x)$ denotes the \emph{fractional anisotropy index} \cite{Hagmann2006} given by
\begin{linenomath}\begin{equation*}
	FA(x) := \sqrt{\dfrac{(\lambda_{1} - \lambda_{2})^2 + (\lambda_{2} - \lambda_{3})^2 + (\lambda_{1} - \lambda_{3})^2}{2 (\lambda_{1}^2 + \lambda_{2}^2 + \lambda_{3}^2)}}.
\end{equation*}\end{linenomath}
Finally, it is well known that the malignant glioma cells positioned in the white matter exercise quicker motility than those situated in the grey matter \cite{Chicoine1995,Kelly1994,Silbergeld1997,Swanson2000}. To account for this effect, in \eqref{Tumour_Diffusion_Tensor}, we use a regulator term $D_{WG}(\cdot)$ that is given by
\begin{linenomath}
\begin{equation}
	D_{WG}(x) = a+(1-a)\Big(\big( D_{G}\: g(x) + w(x)\big) \ast \psi_{\rho}\Big)(x), \label{Eq:Dwg}
\end{equation}
\end{linenomath}
where $0 \leq D_{G} \leq 1$ is the grey matter regulator coefficient, $\ast$ is the convolution operator \cite{Damelin_2011}, $\psi_{\rho}:=\psi(x)/ \rho^{\re{N}}$ denotes the standard mollifier and $g(x)$ and $w(x)$ are the grey and white matter densities provided by the \ds{T1 weighted image} (following an image segmentation process). \re{Finally, $0\leq a\leq 1$ is a model switching parameter distinguishes between different cases (see Section~\ref{sec:Results}).}

In addition, the movement of the cancer cells \re{is} further biased by various adhesion mediated process \cite{Chen_2011,Condeelis2006,Huda2018,Petrie2009,Weiger2013,Wu3949}. Due to the increasing evidence that gliomas induce a fibrous environment within the brain \cite{Gondi2004,Gregorio2018,Kalinin2020,Mohanam1999,Persson2015,Pointer2016,Pullen2018,Ramachandran2017,Veeravalli2012,Veeravalli2012a,Young2009}, in \eqref{Cancer_Cell_Dynamics} we model the overall adhesion process using a non-local flux term that was introduced in \cite{Shuttleworth_2019} (see also \cite{Armstrong2006,Domschke_et_al_2014,Gerisch2008,Peng2016,Shuttleworth2020a,Shuttleworth2020b,Suveges_2020,Suveges_2021} \re{for similar terms}). Specifically, we explore the adhesive interactions of the cancer cells at $x \in \Omega(t)$ with other \re{neighbouring} cancer cells, with the distribution of the non-fibre ECM phase \cite{Ghosh2017,Gras2009,Gras2008,Jacob2016} as well as with the oriented fibre ECM phase \cite{Wolf2009,Wolf_Friedl_2011}, \re{all located} within a sensing region $\Bila(x, R)$ of radius $R > 0$. For this, we define the non-local flux term as
\begin{linenomath}\begin{equation}
	\begin{split}
		\A_{c}(x,t,\bu,\theta_{f}) := \frac{1}{R} \int\limits_{\Bila(0,R)} \K(y) \Big [ & n(y) \big( \S_{cc} c(x+y,t) + \S_{cl} l(x+y,t) \big) \\
		& + \widehat{n}(y, \theta_{f}(x+y,t)) \S_{cF} F(x+y,t) \Big] \big[ 1 - \rho(\bu) \big]^{+} \; dy,
	\end{split}
	\label{Cancer_Adhesion}
\end{equation}\end{linenomath}
where $\S_{cc}$, $\S_{cl}$, $\S_{cF} > 0$ are the strengths of the cell-cell, cell-non-fibre ECM and cell-fibre ECM adhesions, respectively. While we take both $\S_{cl}$ and $\S_{cF}$ as positive constants, we consider the emergence of strong and stable cell-cell adhesion bonds to be positively correlated with the level of extracellular $Ca^{+2}$ ions (one of the non-fibre ECM component) \cite{Gu2014,Hofer2000}. Hence, following \re{the approach in} \cite{Shuttleworth_2019,Shuttleworth2020b,Shuttleworth2020a,Suveges_2020,Suveges_2021}, we describe the cell-cell adhesion strength by
\begin{linenomath}\begin{equation*}
	\S_{cc}(x,t) :=\S_{min} + (\S_{max} - \S_{min}) \exp \bigg[ 1-\dfrac{1}{1-(1-l(x, t))^{2}} \bigg], 
\end{equation*}\end{linenomath}
where $S_{min} > 0$ and $S_{max} > 0$ are the minimum and maximum levels of $Ca^{+2}$ ions. Furthermore, in \eqref{Cancer_Adhesion} we denote by $n(\cdot)$ and $\widehat{n}(\cdot, \cdot)$ the unit radial vector and the unit radial vector biased by the oriented ECM fibres \cite{Shuttleworth_2019,Shuttleworth2020b,Shuttleworth2020a,Suveges_2020,Suveges_2021} defined by
\begin{linenomath}\begin{equation*}
	\begin{split}
		n(y) := &
		\begin{cases}
			\dfrac{y}{\nor{y}_{2}} & \text{if } y \in \Bila(0, R) \setminus \{0\}, \\
			\hfil 0 & \text{if } y = 0,
		\end{cases}
		\\
		\widehat{n}(y, \theta_{f}(x+y)) := &
		\begin{cases}
			\dfrac{y + \theta_{f}(x+y,t)}{\nor{y + \theta_{f}(x+y,t)}_{2}} & \text{if } y \in \Bila(0, R) \setminus \{0\}, \\
			\hfil 0 & \text{if } y = 0,
		\end{cases}
	\end{split}
\end{equation*}\end{linenomath}
respectively (for details on the fibre orientation $\theta_{f}$ see Section~\ref{sec:Fibre_Micro_Scale}). Also, to account for the gradual weakening of all adhesion bonds as we move away from the centre point $x$ within the sensing region $\Bila(x,t)$ in \eqref{Cancer_Adhesion}, we use a radially symmetric kernel $\K(\cdot)$ \cite{Suveges_2020,Suveges_2021} given by
\begin{linenomath}\begin{equation*}
	\K(y) = \psi \Big( \frac{y}{R} \Big), \qquad \forall y\in \Bila(0,R),
\end{equation*}\end{linenomath}
where $\psi(\cdot)$ is the standard mollifier. Finally, in \eqref{Cancer_Adhesion} a limiting term $[1 - \rho(\bu)]^{+} := \max(0, 1 - \rho(\bu))$ is used to prevent the contribution of overcrowded regions to cell migration \cite{Gerisch2008}. \sz{For a \khi{schematic} of this adhesion process, we refer the reader to Figure~\ref{fig:Model_Schematics}(a).}

\subsubsection{Two phase ECM macro-scale dynamics}
In addition to the cancer cell population, the rest of the macro-scale tumour dynamics \khi{are} described by the two-phase ECM. Here, both fibres and non-fibres ECM phases are assumed to be simply described by a degradation term due to the cancer cell population. Hence, per unit time, their dynamics is governed by
\begin{linenomath}\begin{equation}
	\begin{split}
		\dfrac{\partial F}{\partial t} = & - \beta_{F} cF, \\
		\dfrac{\partial l}{\partial t} = & - \beta_{l} cl,
	\end{split}
	\label{Both_ECM_Dynamics}
\end{equation}\end{linenomath}
where $\beta_{F} > 0$ and $\beta_{l} > 0$ are the degradation rates of the fibre and non-fibre ECM phases, respectively.

\subsubsection{The complete macro-dynamics}
\re{In summary, equations \eqref{Cancer_Cell_Dynamics} for cancer cells dynamics and \eqref{Both_ECM_Dynamics} for the two-phase ECM dynamics lead to the following non-dimensional PDE system describing the evolution of tumour at macro-scale:}
\begin{linenomath}\begin{equation}
	\begin{split}
		\dfrac{\partial c}{\partial t} = & \nabla \nabla :\big[ \DD_{T}(x)c \big] - \nabla \cdot \big[ c \A_{c}(x,t,\bu,\theta_{f}) \big] + \mu c \big[ 1 - \rho(\bu) \big], \\
		\dfrac{\partial F}{\partial t} = & - \beta_{F} cF, \\
		\dfrac{\partial l}{\partial t} = & - \beta_{l} cl.
	\end{split}
	\label{Full_Macro_Scale_Dynamics}
\end{equation}\end{linenomath}
\re{To complete the macro-scale model description, we consider} zero-flux boundary conditions and appropriate initial conditions (for instance the ones given in Section~\ref{sec:Results}).

\subsection{Micro-Scale Processes and the Double Feedback Loop}
Since the cancer invasion process is genuinely a multi-scale phenomenon, several micro-scale processes are closely linked to the macro-scale dynamics \cite{Weinberg2006}. In this work, we consider two of these micro-processes, namely the rearrangement of the ECM fibres micro-constituents \cite{Shuttleworth_2019} and the cell-scale proteolytic processes that occur at the leading edge of the tumour \cite{Dumitru_et_al_2013}. Here we briefly outline these micro-processes, in addition to the naturally arising double feedback loop that ultimately connects the micro-scale and the macro-scale.

\subsubsection{Two-scale representation and dynamics of fibres}
\label{sec:Fibre_Micro_Scale}
To represent the oriented fibres on the macro-scale, we follow \cite{Shuttleworth_2019}. \re{There}, the authors characterised not only the amount of fibres $F(x,t)$, but also their ability to withstand incoming cell fluxes and forces through their spatial bias. By considering a cell-scale micro-domain $\delta Y(x) := x + \delta Y$ of appropriate micro-scale size $\delta > 0$, both of these characteristics are induced by the microscopic fibre distribution $f(z,t)$, with $z \in \delta Y(x)$. In fact, both of them are captured though a vector field representation $\theta_{f}(x,t)$ of the ECM micro-fibres \cite{Shuttleworth_2019} that is defined as:
\begin{linenomath}\begin{equation}
	\theta_{f}(x,t) := \dfrac{1}{\lambda (\delta Y(x))} \int\limits_{\delta Y(x)} f(z,t) dz \cdot \dfrac{\theta_{f,\delta Y(x)}(x,t)}{\nor{\theta_{f,\delta Y(x)}(x,t)}_{2}},
	\label{Fibre_Orientation}
\end{equation}\end{linenomath}
where $\lambda(\cdot)$ is the Lebesgue measure in $\R^{d}$ and $\theta_{f,\delta Y(x)}(\cdot,\cdot)$ is the revolving barycentral orientation given by \cite{Shuttleworth_2019}
\begin{linenomath}\begin{equation*}
	\theta_{f,\delta Y(x)}(x,t) := \dfrac{\int\limits_{\delta Y(x)} f(z,t) (z-x) dz}{\int\limits_{\delta Y(x)} f(z,t) dz}.
\end{equation*}\end{linenomath}
Hence, the fibres' ability to withstand forces is naturally defined by this vector field representation \eqref{Fibre_Orientation} and their amount distributed at a macro-scale point $(x,t)$ is given by
\begin{linenomath}\begin{equation*}
	F(x,t) := \nor{\theta_{f}(x,t)}_{2},
\end{equation*}\end{linenomath}
which is precisely the mean-value of the micro-fibres distributed on $\delta Y(x)$. Since both of these macro-scale oriented ECM fibre characteristics ($F(x,t)$ and $\theta_{f}(x,t)$) that we use in the macro-scale dynamics \eqref{Full_Macro_Scale_Dynamics}, genuinely emerge from the micro-scale distribution of the ECM fibres $f(z,t)$, we refer to this link as the \emph{fibres bottom-up} link. \sz{An illustration of a micro-domain $\delta Y(x)$ and its corresponding macro-scale orientation $\theta_{f}(x,t)$ can be seen in Figure~\ref{fig:Model_Schematics}(b).}

On the other hand, there is also a naturally arising link that connects the macro-scale to this micro-scale, namely the \emph{fibres top-down} link. This connection is initiated by the movement of the cancer cell population that trigger a rearrangement of ECM fibres micro-constituents on each micro-domain $\delta Y(x)$ (enabled by the secretion of matrix-degrading enzymes that can break down various ECM proteins). Hence, using the fact that the fully anisotropic diffusion term can be rewritten as $\nabla \nabla :\big[ \DD_{T}(x)c \big] = \nabla \cdot [\DD_{T}(x) \nabla c + c \nabla \cdot \DD_{T}(x)]$, the fibre rearrangement process is kicked off by the cancer cell spatial flux
\begin{linenomath}\begin{equation}
	\F_{c}(x,t) := \DD_{T}(x) \nabla c + c \nabla \cdot \DD_{T}(x) - c\A_{c}(x,t,\bu,\theta_{f}),
	\label{Cancer_Cell_Spatial_Flux}
\end{equation}\end{linenomath}
which is generated by the tumour macro-dynamics \eqref{Full_Macro_Scale_Dynamics}. Then, at any spatio-temporal point $(x,t) \in \Omega(t) \times [0,T]$ this flux \eqref{Cancer_Cell_Spatial_Flux} gets naturally balanced in a weighted fashion by the macro-scale ECM fibre orientation $\theta_{f}(\cdot, \cdot)$, resulting in a \emph{rearrangement flux} \cite{Shuttleworth_2019}
\begin{linenomath}\begin{equation}
	r(\delta Y(x),t) := \omega(x,t)\F_{c}(x,t) + (1-\omega(x,t))\theta_{f}(x,t),
	\label{Fibre_Rearrangement_Vector}
\end{equation}\end{linenomath}
with $\omega(x,t) := c(x,t) / (c(x,t) + F(x,t))$, that acts uniformly upon the micro-fibre distribution on each micro-domain $\delta Y(x)$. Ultimately, this macro-scale rearrangement vector \eqref{Fibre_Rearrangement_Vector} induces a micro-scale reallocation vector $\nu_{\delta Y(x)}(z,t)$ \cite{Shuttleworth_2019}, enabling us to appropriately calculate the new position $z^{*}$ of any micro-node $z \in \delta Y(x)$ as
\begin{linenomath}\begin{equation}
	z^{*} := z + \nu_{\delta Y(x)}(z,t).
	\label{Fibre_Rearrangement_New_Position}
\end{equation}\end{linenomath}
For further details on the micro-fibre rearrangement process, we refer the reader to Appendix~\ref{sec:Appendix_Further_Details_On_The_Fibre_Rearrangement} and \cite{Shuttleworth_2019,Shuttleworth2020a,Shuttleworth2020b,Suveges_2020,Suveges_2021}.

\subsubsection{MDE micro-dynamics and its links}
The second micro-scale process that we take into consideration is the proteolytic molecular process that occurs along the invasive edge of the tumour and is driven by the cancer cells' ability to secrete several types of \emph{matrix-degrading enzymes} (MDEs) (for instance, matrix-metalloproteinases) within the proliferating rim \cite{Hanahan2000,Hanahan2011,Lu2011,Parsons1997,Pickup2014}. Subsequent to the secretion, these MDEs are subject to spatial transport within a cell-scale neighbourhood of the tumour interface and, as a consequence, they degrade the peritumoral ECM, resulting in changes of tumour boundary morphology \cite{Weinberg2006}.

To explore such \khi{a} micro-scale process, we adopt the approach that was first introduced in \cite{Dumitru_et_al_2013} where the emerging spatio-temporal MDEs micro-dynamics is considered on an appropriate cell-scale neighbourhood of the tumour boundary $\partial \Omega(t)$. This neighborhood is represented by a time-dependent bundle of overlapping cubic micro-domains $\{ \epsilon Y \}_{\epsilon Y \in \P(t)}$, with $\epsilon > 0$ being the size of each micro-domain $\epsilon Y$, which allows us to decompose the overall MDE micro-process, transpiring on $\bigcup_{\epsilon Y \in \P(t)} \epsilon Y$, into a union of proteolytic micro-dynamics occurring on each $\epsilon Y$; \sz{see also Figure~\ref{fig:Model_Schematics}(c)}. Hence, choosing an arbitrary micro-domain $\epsilon Y$ and a macroscopic time instance $t_{0} \in [0,T]$, we follow the evolution of the MDE micro-dynamics during the time period $[t_{0}, t_{0} + \Delta t]$, with appropriately chosen $\Delta t > 0$ and within the associated micro-domain $\epsilon Y$. By denoting the spatio-temporal distribution of the MDEs by $m(y,\tau)$ at any micro-point $(y,\tau) \in \epsilon Y \times [0, \Delta t]$, we observe that the cancer cell population, located within an appropriately chosen distance $\gamma_{h} > 0$ from $y \in \epsilon Y$, induce a source $h(y, \tau)$ of MDEs which can be mathematically described via a non-local term \cite{Dumitru_et_al_2013}
\begin{linenomath}\begin{equation}
	h(y,\tau) = 
	\begin{cases}
		\dfrac{\int\limits_{\Bila(y, \gamma_{h}) \cap \Omega (t_{0})} c(x,t_{0} + \tau) \; dy}{\lambda (\Bila(y, \gamma_{h}) \cap \Omega (t_{0}))} & y \in \epsilon Y \cap \Omega(t_{0}), \\[10pt]
		\hfill 0 & y \notin \epsilon Y \setminus (\Omega(t_{0}) + \{ z \in Y \; | \; \| z \|_{2} < \rho \}),
	\end{cases}
	\label{MDE_Source}
\end{equation}\end{linenomath}
where $0 < \rho < \gamma_{h}$ is a small mollification range and $\Bila(y, \gamma_{h})$ denotes the $\nor{\cdot}_{_{\infty}}$ ball of radius $\gamma_{h}$ centred at a micro-node $y$. Since the calculation of this micro-scale MDE source \eqref{MDE_Source} directly involves the macro-scale cancer cell population $c(\cdot, \cdot)$, we observe a naturally arising \emph{MDE top-down} link that connects the macro-scale to the MDE micro-scale. In fact, such source term \eqref{MDE_Source} allows us to describe the spatio-temporal evolution of the MDEs micro-scale distribution $m(\cdot, \cdot)$ by \cite{Dumitru_et_al_2013}
\begin{linenomath}\begin{equation}
    \begin{split}
    	\dfrac{\partial m}{\partial \tau} & = D_{m} \Delta m + h(y, \tau), \\
    	m(y,0) & = 0, \\
    	\dfrac{\partial m}{\partial n} \Big |_{\partial \epsilon Y} & = 0,
    \end{split}
	\label{MDE_Equation}
\end{equation}\end{linenomath}
where $D_{m} > 0$ is the constant MDEs diffusion coefficient and $n$ denotes the outward normal vector. As it was shown in \cite{Dumitru_et_al_2013}, we can use the solution of the MDEs micro-dynamics \eqref{MDE_Equation} to acquire both movement direction and magnitude of a tumour boundary point $x^{*}_{\epsilon Y}$ within the peritumoral area covered by the associated boundary micro-domain $\epsilon Y$. This ultimately causes a boundary movement, and as a consequence we obtain a new evolved tumour macro-domain $\Omega(t_{0} + \Delta t)$, the link of which we refer to as the \emph{MDE bottom-up link}. For \sz{illustration of the boundary movement we refer the reader to Figure~\ref{fig:Model_Schematics}(c) and for} further details of the MDE micro dynamics see Appendix~\ref{sec:Appendix_Further_Details_On_The_MDE_Micro_Scale} or \cite{Dumitru_et_al_2013,Shuttleworth_2019,Shuttleworth2020a,Shuttleworth2020b,Suveges_2020,Suveges_2021}.

\section{Computational Results: Numerical Simulations in 3D}\label{sec:Results}
We start this section by briefly discussing the numerical method that we use to solve the macro-scale dynamics \eqref{Full_Macro_Scale_Dynamics}, and for details on the numerical approach used for the two micro-scale dynamics (fibres and MDE) we refer the reader to \cite{Suveges_2020,Suveges_2021}. Here, we use the method of lines approach to discretise the macro-scale tumour dynamics \eqref{Full_Macro_Scale_Dynamics} first in space, and then, for the resulting system of ODEs, we employ a non-local predictor-corrector scheme \cite{Shuttleworth_2019}. In this context, we carry out the spatial discretisation on a uniform grid, where both spatial operators (fully anisotropic diffusion and adhesion) are accurately approximated in a convolution-driven fashion. First, we note that the fully anisotropic diffusion term can be split into two parts
\begin{linenomath}\begin{equation}
    \nabla \nabla :\big[ \DD_{T}(x)c \big] = \underbrace{\nabla \cdot \big[ \DD_{T}(x) \nabla c \big]}_{\text{diffusive}} + \underbrace{\nabla \cdot \big[c  \nabla \cdot \DD_{T}(x) \big]}_{\text{advective}},
    \label{Diffusion_Split_Numerics}
\end{equation}\end{linenomath}
which enables us to use a combination of two appropriate distinct schemes for an accurate approximation. While for the diffusive part in \eqref{Diffusion_Split_Numerics}, we use the \emph{symmetric finite difference} scheme \cite{Es2014,Guenter2005}, for the combination of the advective \eqref{Diffusion_Split_Numerics} and adhesion operators \eqref{Cancer_Adhesion} (\re{i.e.,} $\nabla \cdot \big[ c \big(\A_{c}(x,t,\bu,\theta_{f}) + \nabla \cdot \DD_{T}(x) \big) \big]$) we use the standard \emph{first-order upwind} finite difference scheme which ensures positivity and helps avoiding spurious oscillations in the solution. Finally, to approximate the adhesion integral $\A_{c}(x,t,\bu, \theta_{f})$, \re{we consider an approach similar to \cite{Suveges_2020,Suveges_2021}, and} use $N_{s}$ random points located within the sensing region $\Bila(0,R)$ and sums of discrete-convolutions. 

\subsection{Initial Conditions}
\begin{figure}[t!]
\centering
\includegraphics[width=10.5 cm]{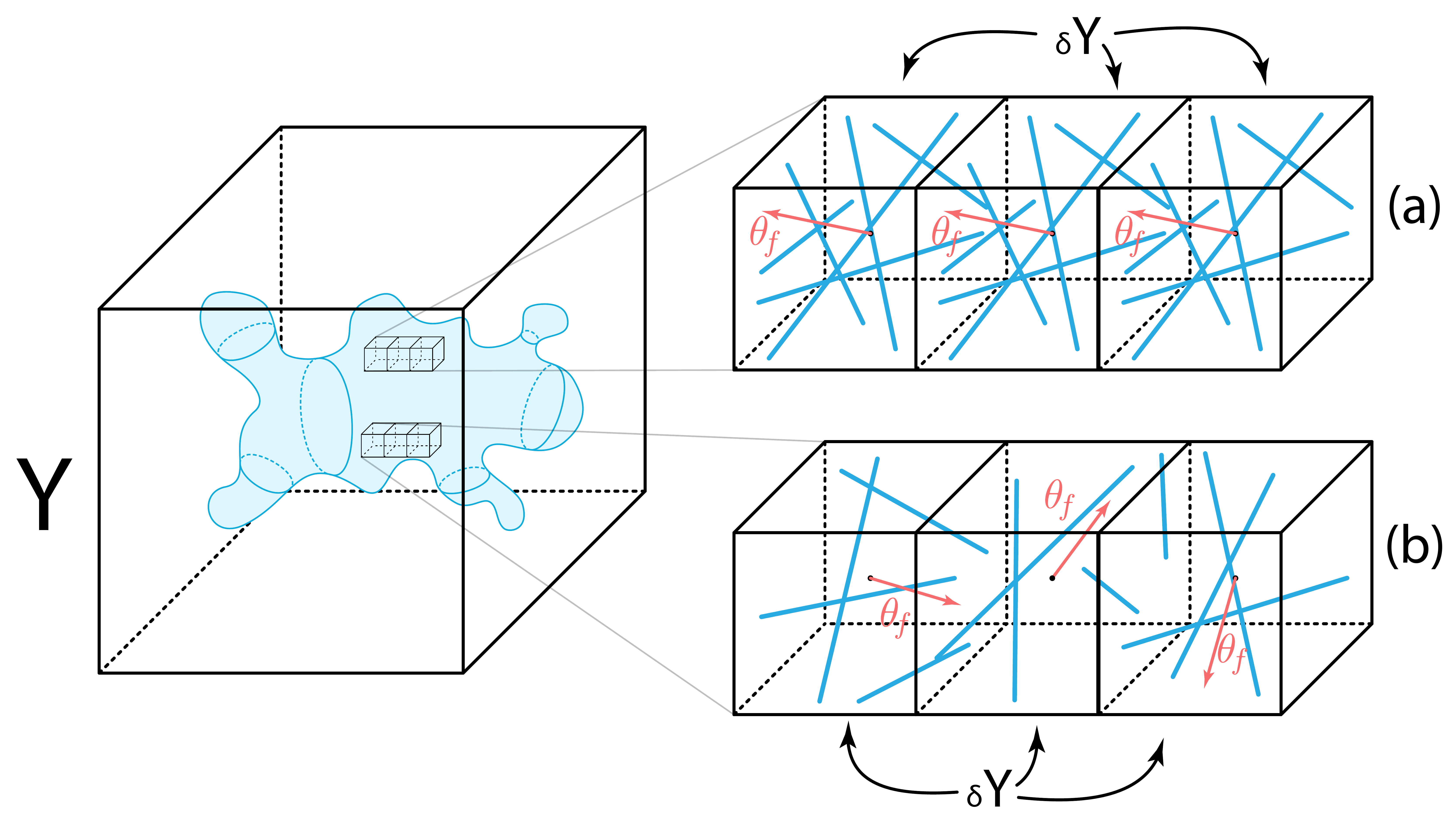}
\caption{\sz{Schematics of the initial condition of the micro-fibres (blue lines) within a micro-domain $\delta Y(x)$ of orientation $\theta_{f}(x,t)$ located in the (a) white matter and in the (b) grey matter.}\label{fig:White_Grey_Fibre_Structure}}
\end{figure}

For the numerical simulations presented in this paper, we consider the tissue cube $Y=[0,4] \times [0,4] \times [0,4]$ with the following initial condition for the cancer cells
\begin{linenomath}\begin{equation*}
	c(x,0) = \dfrac{1}{2} \exp\bigg( \dfrac{-\nor{x}_{2}^{2}}{0.02} \bigg) \cdot \chi_{_{\Bila((2,2,2),0.25)}},
\end{equation*}\end{linenomath}
and for the non-fibre ECM phase, the initial condition $l(x,0)$ is acquired by appropriately scaling the \ds{T1 weighted image} \re{via a normalising constant}. Current DTI \ds{scans} do not provide suitable resolution to determine the underlying micro-fibre distributions, and so here, we describe the initial micro-fibre distribution within a micro-domain $\delta Y(x)$ as follows. When the macro-scale point $x$ that corresponds to the micro-domain $\delta Y(x)$ is located in the grey matter, then within $\delta Y(x)$ we randomly draw straight lines until the ratio between \re{the} points that belong and the points that do not belong to the collection of lines is about $35\%$ : $65\%$. On the other hand, when the point $x$ is located within the white matter, we use a set of predefined lines with the same point ratio ($35\%$ : $65\%$), ultimately achieving a random orientation within the grey matter and an aligned orientation within the white matter \ds{\cite{Raffelt2017}}. Finally, the grey matter's fibre density is assumed to be $1/D_{G}$ times smaller than the density in the white matter \ds{\cite{Raffelt2017}}. \sz{A schematics of this initial condition for the micro-fibres can be seen in Figure~\ref{fig:White_Grey_Fibre_Structure}.} Hence, we also incorporate the information about the white and grey matter tracks provided by the \ds{T1 weighted image} into our micro-scale fibre distribution.

\subsection{Numerical Simulations in 3D}
Here, we present the 3D numerical solutions of the multi-scale model described above, \re{for the} parameter values \re{listed in} Table~\ref{tab:Parameter_Set} in Appendix~\ref{sec:Appendix_Parameter_Set} (any alteration from these values will be stated accordingly). To display the advanced tumours at time $50 \Delta t$, we show four panels for each simulation results. \re{In the first three panels we show} the three \re{classical} cross-section \re{planes} \emph{i.e.,} the coronal \re{plane} (the head of the subject is viewed from behind), the axial \re{plane} (the head of the subject is viewed from above) and the sagittal \re{plane} (the head of the subject is viewed from the left). \re{In the last} panel of each simulation we show the 3D image of the brain with the embedded tumour alongside the 3D tumour in isolation.

The three Figures shown below investigate tumour evolution when the initial tumour starts in different regions of the brain.

\begin{figure}[t!]
\widefigure
\includegraphics[width=15 cm]{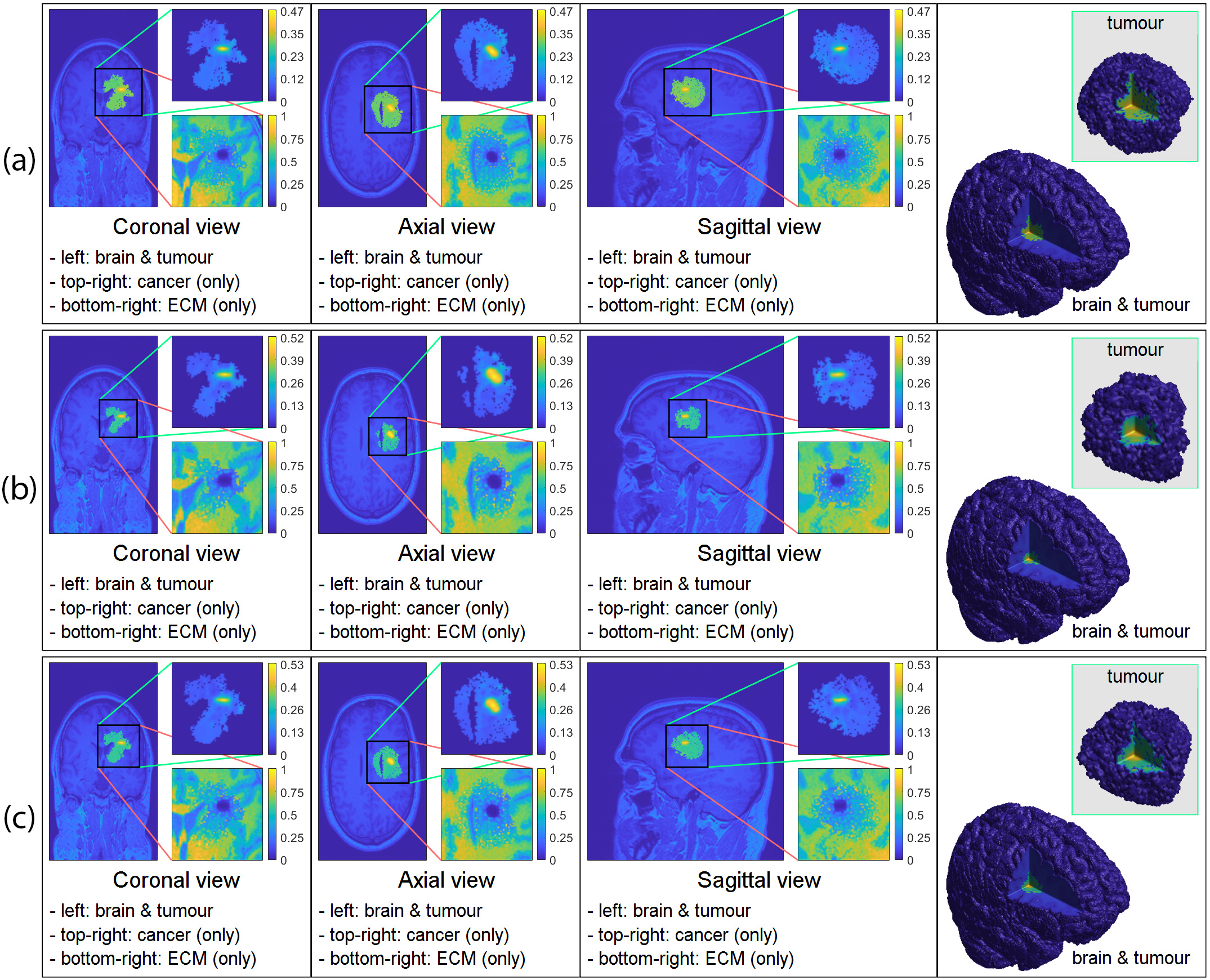}
\caption{3D computer simulation results (a) with only white-grey matter dependency ($r=1$), (b) with only DTI data used ($D_{G} = 1$) (c) with both white-grey matter dependency and DTI data incorporated. To present the simulations, we divide each result into four panels: coronal, axial, sagittal and 3D view. Within each coronal, axial and sagittal views, we show the tumour embedded within the brain on the left, the cancer cell density on the top-right and the ECM density on the bottom-right. In the 3D view (the most \khi{right} panel in each results) we show the cross-section of the whole brain with the tumour on the bottom-left corner and on the top-right \khi{corner} we show the isolated tumour. \label{fig:Simulations_1}}
\end{figure}

\begin{figure}[t!]
\widefigure
\includegraphics[width=15 cm]{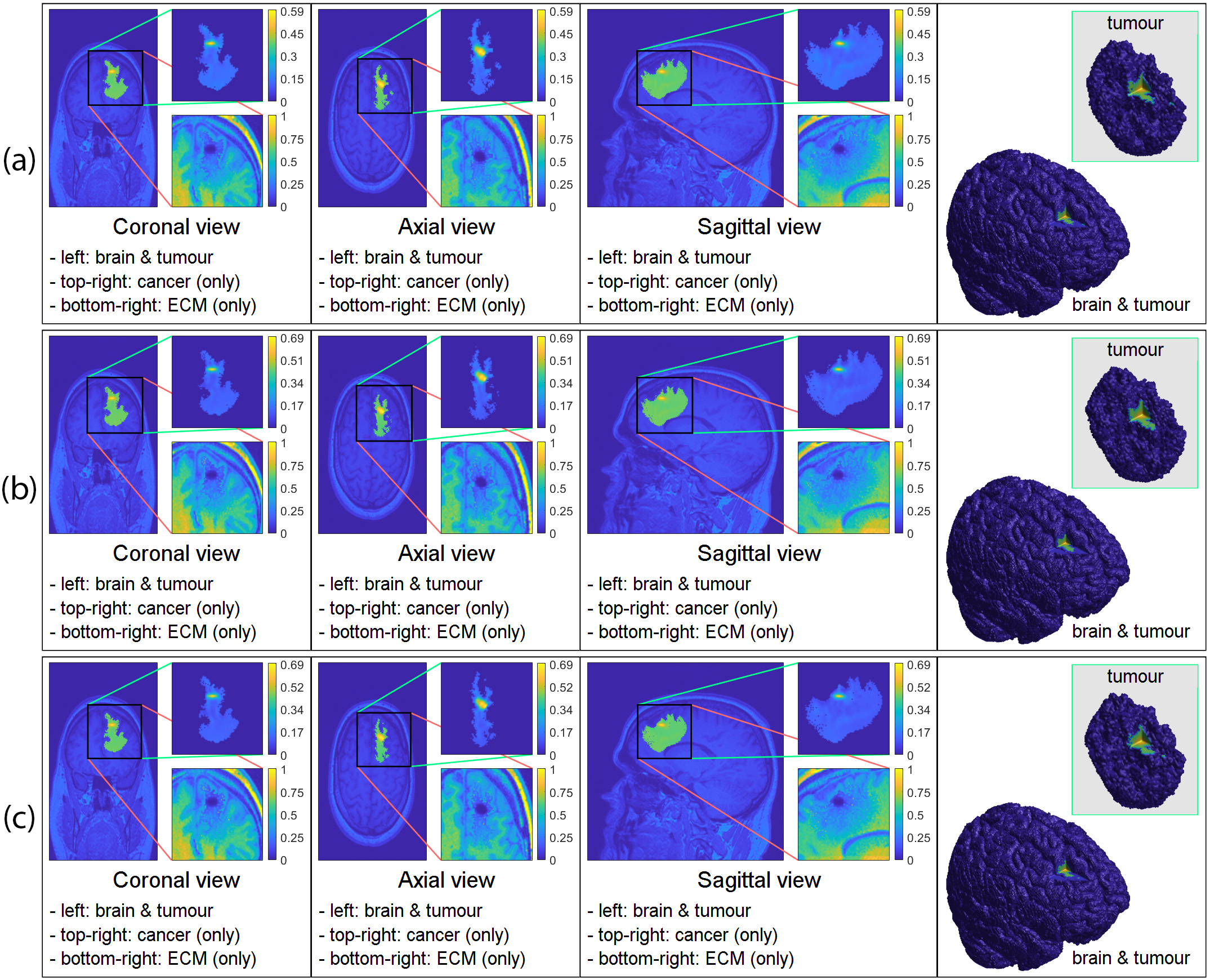}
\caption{3D computer simulation results (a) with only white-grey matter dependency ($r=1$), (b) with only DTI data used ($D_{G} = 1$) (c) with both white-grey matter dependency and DTI data incorporated. To present the simulations, we divide each result into four panels: coronal, axial, sagittal and 3D view. Within each coronal, axial and sagittal views, we show the tumour embedded within the brain on the left, the cancer cell density on the top-right and the ECM density on the bottom-right. In the 3D view (the most \khi{right} panel in each results) we show the cross-section of the whole brain with the tumour on the bottom-left corner and on the top-right \khi{corner} we show the isolated tumour. \label{fig:Simulations_2}}
\end{figure}

\begin{figure}[t!]
\widefigure
\includegraphics[width=15 cm]{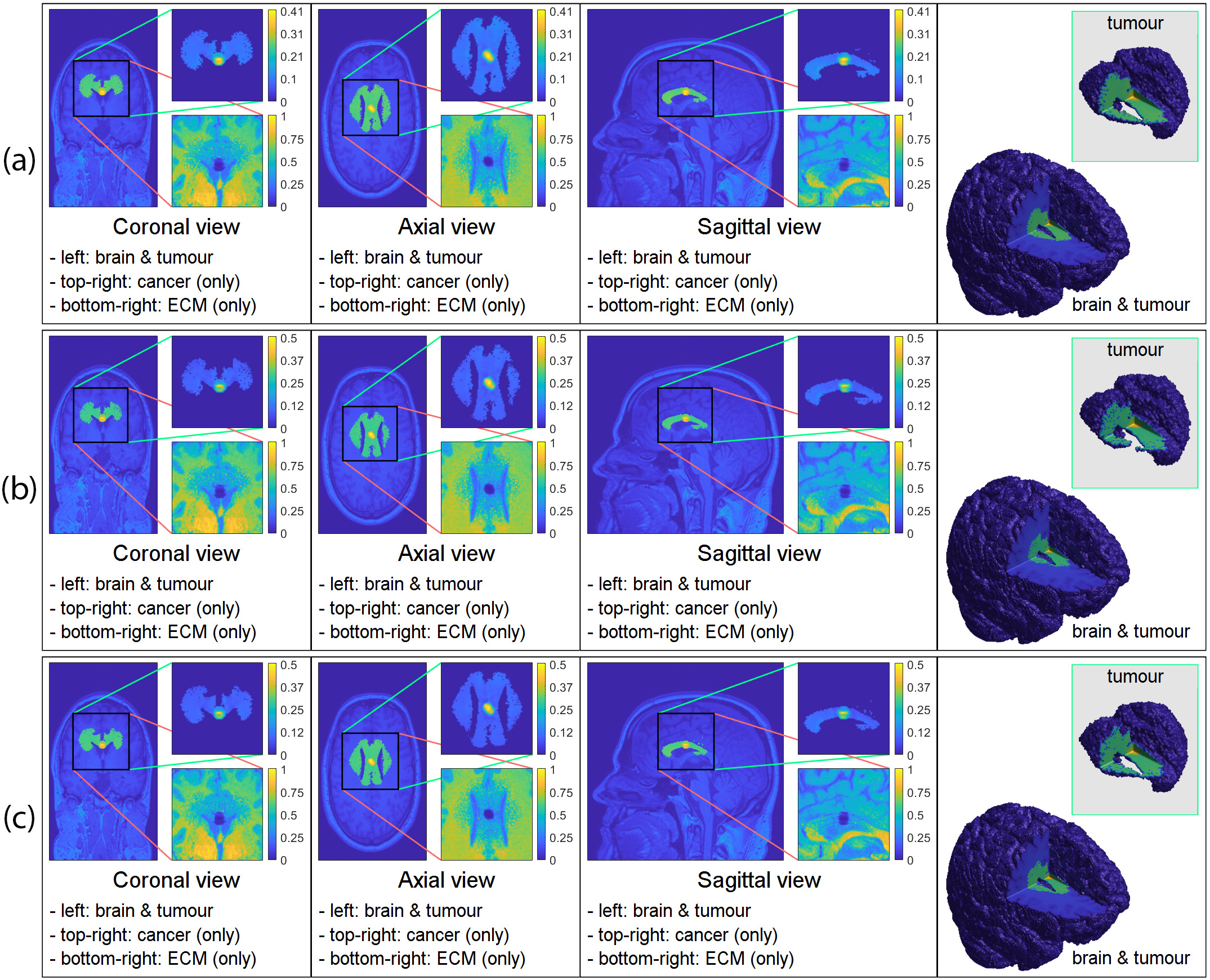}
\caption{3D computer simulation results (a) with only white-grey matter dependency ($r=1$), (b) with only DTI data used ($D_{G} = 1$) (c) with both white-grey matter dependency and DTI data incorporated. To present the simulations, we divide each result into four panels: coronal, axial, sagittal and 3D view. Within each coronal, axial and sagittal views, we show the tumour embedded within the brain on the left, the cancer cell density on the top-right and the ECM density on the bottom-right. In the 3D view (the most \khi{right} panel in each results) we show the cross-section of the whole brain with the tumour on the bottom-left corner and on the top-right \khi{corner} we show the isolated tumour. \label{fig:Simulations_3}}
\end{figure}

In Figure~\ref{fig:Simulations_1} we present three distinct cases obtained by varying different parameters that appear in the tumour diffusion tensor $\DD_{T}(x)$ defined in \eqref{Tumour_Diffusion_Tensor}. In Figure~\ref{fig:Simulations_1} (a) \re{we assume that the tensor $\DD_{T}(x)$ depends on the white-grey matter and for that purpose we set $r=1$ in \eqref{Tumour_Diffusion_Tensor} and $a=0$ in \eqref{Eq:Dwg};} this results in isotropic tumour diffusion. In Figure~\ref{fig:Simulations_1} (b) we use the DTI data (\re{i.e., there is no a-priory assumption about the preferential direction for cell movement in white matter}) and thus we set $a=1$ in \eqref{Eq:Dwg} \re{(with $r=0.1$, as in Table~\ref{tab:Parameter_Set})}; this results in an anisotropic diffusion that does not depend explicitly on the white-grey matter. In Figure~\ref{fig:Simulations_1} (c) we use both DTI data and the white-grey matter dependency (i.e., $r=0.1$ and $a=0$), with the baseline parameters from Table~\ref{tab:Parameter_Set}. Here, it is worth mentioning that even though we do not use the \ds{T1 weighted image} \re{to obtain functions $w(x)$ and $g(x)$ that appear in $D_{WG}$ (as $D_{WG}=1$ in Figure~\ref{fig:Simulations_1} (b), since $a=1$)} we still use the \ds{T1 weighted image} to initialise the micro-scale non-fibre initial density as well as the initial micro-scale fibre distributions as described above.\\ 
In \re{all} these simulations \re{shown in Figure~\ref{fig:Simulations_1}}, we place the small initial tumour in the middle-right part of the brain, and we show the results of the three cases at time $50 \Delta t$ where we observe significant tumour morphology changes across the three cases. \re{By} comparing Figure~\ref{fig:Simulations_1} (a) to (b) and Figure~\ref{fig:Simulations_1} (b) to (c), \re{we see that} when we include the white-grey matter dependency function $D_{WG}$ within the tumour diffusion tensor it leads to a more advanced tumour. On the other hand, comparing Figure~\ref{fig:Simulations_1} (a) to (c) shows that including the DTI data, which creates an anisotropic tumour diffusion term, leads to a slight reduction in tumour spread. Furthermore, in all three cases, we can notice that the advancing tumour tends to mostly follow the white matter tracks and usually avoids the invasion of tissues located in the grey matter. This invasion resulted in the degradation (and rearrangement) of the ECM that we can see in the bottom-right of each panel (coronal, axial and sagittal) which enabled the tumour to further expand into the surrounding tissues.
\end{paracol}
\begin{paracol}{2}
\linenumbers
\switchcolumn

In Figure~\ref{fig:Simulations_2} we keep the same three cases as in Figure~\ref{fig:Simulations_1}, \emph{i.e.,} Figure~\ref{fig:Simulations_2} (a) only white-grey matter dependency, Figure~\ref{fig:Simulations_2} (b) only DTI data and Figure~\ref{fig:Simulations_2} (c) both. However, here we place the initial tumour in the front-right part of the brain and show the results of tumour invasion at the final time $50 \Delta t$. Due to the initial position of the tumour, we can see a tumour that \khi{is growing} away from the skull towards the centre of the brain as well as it \khi{is} mainly \khi{following} the white matter. This creates a highly heterogeneous elongated tumour with many branching outgrowths. On the other hand, in Figure~\ref{fig:Simulations_2} we only see slight differences between the three cases. This contradicts the results from Figure~\ref{fig:Simulations_1} and suggests that both the DTI data and white-grey matter dependency may not always be decisive factor in tumour morphology.
\end{paracol}
\begin{paracol}{2}
\linenumbers
\switchcolumn

Similarly to Figure~\ref{fig:Simulations_1} and Figure~\ref{fig:Simulations_2}, in Figure~\ref{fig:Simulations_3} we keep the same three cases (Figure~\ref{fig:Simulations_3} (a) only white-grey matter dependency, Figure~\ref{fig:Simulations_3} (b) only DTI data and Figure~\ref{fig:Simulations_3} (c) both) while we place the initial tumour mass in the middle of the brain and present the results at time $50 \Delta t$. As a consequence of the initial location, we see a "butterfly" shaped tumour that branched to both the left and right side of the brain with some asymmetry. Also, as in Figure~\ref{fig:Simulations_2} we can see that all three cases are quite similar, and so the additional information provided by both the DTI data and white-grey matter dependency seems to be unnecessary for this initial condition. However, we must note that the initial conditions (fibre and non-fibre ECM) still uses the information provided by the \ds{T1 weighted image}, and so here, we only investigate the effect of changing the diffusion tensor. 

\end{paracol}
\begin{paracol}{2}
\linenumbers
\switchcolumn

As we mentioned, we see significant differences between the three cases only in Figure~\ref{fig:Simulations_1}. This either indicates that the anisotropic diffusion tensor provides valuable information only in certain cases or that the initial micro-fibre density differs from the one that produced the DTI \ds{scan} (\emph{i.e.,} the actual distribution). Since we use an artificial micro-fibre structure that does not depend on the DTI \ds{scan} which also aid the movement of the cancer cell population via the adhesion integral $\A_{c}(\cdot, \cdot, \cdot, \cdot)$ defined in \eqref{Cancer_Adhesion}, it is possible that in this specific case the micro-scale fibre distribution introduced a significantly different travelling direction than the DTI data, resulting in discrepancies between the simulations. However, due to the resolution of current DTI \ds{scans}, it is not possible to construct a unique fibre distribution within a micro-domain $\delta Y(x)$. Hence, to genuinely capture the underlying brain structures that we can use within a mathematical model, our results suggest that DTI \ds{scans} with \ds{their} present resolution may not be sufficient, and one might need to look into either obtaining better resolution DTI \ds{scans} or \re{combine this with the strength of} different technologies such as magnetic resonance elastography. Nonetheless, this exceeds the current scope of this work and requires further investigation.

\section{Discussion and Final Remarks}\label{sec:Discussion}
In this study, we have further extended the \re{2D multi-scale moving-boundary framework previously introduced in }\cite{Dumitru_et_al_2013,Shuttleworth_2019}, \re{by developing it to 3D and applying it to the study of glioma invasion within the brain}. Since experiments are limited within the brain, we focused on incorporating \ds{DTI and T1 weighted scans} into our framework \re{to} provide insights into the structure of the brain, the tumour, and the surrounding tissue.

The original framework developed in \cite{Dumitru_et_al_2013,Shuttleworth_2019} modelled a generic tumour in a 2D setting, and so to model gliomas within a 3D brain, we extended this modelling approach by considering the structural information provided by both DTI and \ds{T1 weighted scans}. We used \re{both DTI and \ds{T1 weighted scans}} to construct the tumour diffusion tensor $\DD_{T}(x)$ defined in \eqref{Tumour_Diffusion_Tensor}, \re{which} resulted in a fully anisotropic diffusion term. While the \ds{T1 weighted image} \re{can} give different diffusion rates based on whether the cancer cells are located in the white or grey matter, the DTI data is used to incorporate the underlying brain structure and to give higher diffusion rates along specific directions based on how the measured water molecules behaved within the brain. \re{The \ds{T1 weighted image}, which provided the} white-grey matter densities, were also used in our initial conditions for both ECM phases. Hence, the initial density of the non-fibre ECM phase was taken as a \re{normalised} version of the \ds{T1 weighted image}, and the initial condition of the micro-fibre distribution and magnitude were also considered to be dependent on the white-grey matter structure. \re{Furthermore, as the available DTI \ds{scans} lack the adequate resolution required to construct more appropriate micro-fibre distributions, in this work we considered a simple case where we set the fibre distributions to be either random or oriented based on whether they are positioned in the grey or white matter, respectively.} 

We used this new \re{3D} model to explore the effects of the anisotropic diffusion term for the cancer cell population. Our numerical simulations in Figure~\ref{fig:Simulations_1} \re{showed that} including an anisotropic diffusion term may lead to significant changes in the overall tumour morphology. \re{However, it seems that these changes depend on the position of the tumour inside the brain, as Figures~\ref{fig:Simulations_2} and~\ref{fig:Simulations_3} do not exhibit changes consistent to the ones observed between the three sub-panels of Figure~\ref{fig:Simulations_1}. This may be the result of the underlying brain structure and its microscopic fibre representation, which seems to take a leading role in influencing cancer-invasion patterns through the underlying cell-adhesion process (see Eq. \eqref{Cancer_Adhesion}), overshadowing this way the diffusion process. More precisely, the simplified fibre representation might not be sufficient for Figure~\ref{fig:Simulations_1}, where the initial tumour was positioned in the right-middle part of the brain. However, this fibre representation might be enough for Figure~\ref{fig:Simulations_2} (with tumour positioned in the front-right of the brain) and for Figure~\ref{fig:Simulations_3} (with tumour positioned in the middle of the brain), where we did not see significant morphological differences between the three sub-panels considered in each of these figures.}





To conclude this study, we mention that further investigation is needed to determine whether these changes in tumour invasion patterns are caused by the lack of directional information on the fibre micro-scale level or an anisotropic diffusive cell motility is necessary to better represent the invasion process. A feasible approach would also be to use a new imaging technology called magnetic resonance elastography, but this is beyond the scope of this current work. \khi{Finally, as our simulations are able to reproduce known tumour patterns of growth seen clinically, future experiments will be refined by MRI data collected prospectively from glioma patients and also incorporate the effects of their radiotherapy and chemotherapy treatments.}

\vspace{6pt} 



\authorcontributions{ All authors contributed to this work.}

\funding{This research was funded by EPSRC DTA EP/R513192/1.}

\institutionalreview{``Not applicable''.}



\conflictsofinterest{The authors declare no conflict of interest.} 


\abbreviations{The following abbreviations are used in this manuscript:\\

\noindent 
\begin{tabular}{@{}ll}
MRI & Magnetic Resonance Imaging\\
DTI & Diffusion Tensor Imaging\\
ECM & Extracellular matrix\\
MDE & Matrix degrading enzymes\\
PDE & Partial differential equation
\end{tabular}}

\appendixtitles{yes} 
\appendixstart
\appendix
\section{Parameter Values}
\label{sec:Appendix_Parameter_Set}
In Table~\ref{tab:Parameter_Set}, we summarise the parameter values that were used in the presented numerical simulations.

\begin{table}[!h]
\caption{Parameter set used for the numerical simulations.}
\label{tab:Parameter_Set}
\begin{tabular}{lllc}
\hline\noalign{\smallskip}
\textbf{Variable} & \textbf{Value} & \textbf{Description} & \textbf{Reference} \\
\noalign{\smallskip}\hline\noalign{\smallskip}

$D_{c}$ & $1.25 \times 10^{-4}$ & Diffusion coeff. for the cancer cell population & \cite{Painter2013} \\
$D_{G}$ & $0.25$ & Grey matter regulator coefficient & Estimated \\
$r$ & $0.1$ & Degree of randomised turning & \cite{Painter2013} \\
$a$ & $0$ & Model switching parameter & Estimated \\
$\K_{FA}$ & $100$ & Cell's sensitivity to the directional information & \cite{Painter2013} \\
$\S_{max}$ & $0.5$ & Cell-cell adhesion coeff. & \cite{Shuttleworth_2019} \\
$\S_{min}$ & $0.01$ & Minimum level of cell-cell adhesion & \cite{Suveges_2020} \\
$\S_{cl}$ & $0.01$ & Cell-non-fibre adhesion coeff. & \cite{Shuttleworth_2019} \\
$\S_{cF}$ & $0.3$ & Cell-fibre adhesion coeff. & \cite{Domschke_et_al_2014} \\
$\mu$ & $0.25$ & Proliferation coeff. for cancer cell population & \cite{Domschke_et_al_2014}\\
$\beta_{F}$ & $1.5$ & Degradation coeff. of the fibre ECM & \cite{Suveges_2021} \\
$\beta_{l}$ & $3.0$ & Degradation coeff. of the non-fibre ECM & \cite{Suveges_2021} \\
$\beta$ & $0.8$ & Optimal tissue environment controller & \cite{Dumitru_et_al_2013} \\
$R$ & $0.15$ & Sensing radius & \cite{Shuttleworth_2019} \\
$f_{max}$ & $0.636$ & Maximum of micro-fibre density at any point & \cite{Shuttleworth_2019} \\
$\Delta x$ & $0.03125$ & Macro-scale spatial step-size & \cite{Dumitru_et_al_2013} \\
$\epsilon$ & $0.0625$ & Size of a boundary micro-domain $\epsilon Y(x)$ & \cite{Dumitru_et_al_2013} \\
$\delta$ & $0.03125$ & Size of a fibre micro-domain $\delta Y (x)$ & \cite{Shuttleworth_2019} \\
$N_{s}$ & $450$ & Number of random points used for the & Estimated \\
 & & approximation of the adhesion integral $\A_{c}$ & \\

\noalign{\smallskip}\hline
\end{tabular}
\end{table}

\section{Further Details on the Micro-Fibre Rearrangement Process}
\label{sec:Appendix_Further_Details_On_The_Fibre_Rearrangement}
In Section~\ref{sec:Fibre_Micro_Scale}, we highlighted the fact that the rearrangement of micro-fibre distribution $f(z,t)$ within each $\delta Y(x)$ is initiated by the macro-scale cell fluxes, resulting the redistribution of each micro-fibre pixel $z$ to a new position $z^{*}$. To calculate this new position $z^{*}$, we use the so-called \emph{reallocation vector} $\nu_{\delta Y(x)}(z,t)$ which takes into account the \emph{rearrangement vector} $r(\delta Y(x), t)$, defined in \eqref{Fibre_Rearrangement_Vector}, the degree of alignment between $r(\delta Y(x), t)$ and the \emph{barycentral position vector} $x_{dir} := z-x$ and also incorporates the level of fibres at position $z$. Hence, following \cite{Shuttleworth_2019}, we define it as
\begin{linenomath}\begin{equation*}
	\nu_{\delta Y(x)}(z,t) := \big[ x_{dir}(x) + r(\delta Y(x),t) \big] \cdot \dfrac{f(z,t) [f_{max} - f(z,t)]}{f^{*} + \nor{ r(\delta Y(x),t) - x_{dir}(x)}_{2}} \cdot \chi_{ \{f(z,t) > 0\} },
\end{equation*}\end{linenomath}
where $f_{max} > 0$ is the maximum level of fibres, $f^{*} := f(x,t) / f_{max}$ is the saturation level and $\chi_{ \{f(z,t) > 0\} }$ is the characteristic function of the micro-fibres support. To move the appropriate amount of fibres from position $z$ to the new position $z^{*}$, given in \eqref{Fibre_Rearrangement_New_Position}, we also monitor the available amount of free space at this target position $z^{*}$ via a movement probability $p_{move}$ that we define it as
\begin{linenomath}\begin{equation*}
	p_{move} := \max \Bigg(\dfrac{f_{max} - f(z^{*},t)}{f_{max}}, 0 \Bigg).
\end{equation*}\end{linenomath}
Consequently, we transport $p_{move} \cdot f(z,t)$ amount of fibres to the new position $z^{*}$ and the rest $(1 - p_{move}) \cdot f(z,t)$ remains at the original position $z$.

\section{Further Details on the MDE micro-scale}
\label{sec:Appendix_Further_Details_On_The_MDE_Micro_Scale}
Following \cite{Dumitru_et_al_2013}, we briefly detail here the way the MDE micro-dynamics \eqref{MDE_Equation} determines the macro-boundary of the progressed tumour domain $\Omega(t_{0} + \Delta t)$. To that end, on any arbitrary boundary micro-domain $\epsilon Y \in \P(t_{0})$ we consider an appropriate dyadic cubes decomposition $\{\D_{k}\}_{_{\I}}$, and we denote the barycentre of each $\D_{k}$ by $y_{k}$. Then, a subfamily of small dyadic cubes $\{\D_{k}\}_{_{\J^{*}}}$ is sub-sampled by selecting only those dyadic cubes that are furthest away from the boundary point $x^{*}_{\epsilon Y}$ while being located outside of the tumour domain $\Omega(t_{0})$ and carrying an above average mass of MDEs. This enables us to define the associated direction $\eta_{\epsilon Y}$ and displacement magnitude $\xi_{\epsilon Y}$ of the movement, which are given by
\begin{linenomath}\begin{equation*}
	\begin{split}
		\eta_{\epsilon Y (x^{*}_{\epsilon Y})} & := x^{*}_{\epsilon Y} + \nu \sum_{l \in \mathcal{J}^{*}} \Bigg ( \int\limits_{\mathcal{D}_{l}} m(y, \tau) \; dy \Bigg ) \Big ( y_{l} - x^{*}_{\epsilon Y} \Big ), \qquad \nu \in  [0, \infty), \\
		\xi_{\epsilon Y (x^{*}_{\epsilon Y})} & := \sum_{l \in \mathcal{J}^{*}} \dfrac{\int\limits_{\mathcal{D}_{l}} m(y, \tau) \; dy}{\sum\limits_{l \in \mathcal{J}^{*}} \int\limits_{\mathcal{D}_{l}} m(y, \tau) \; dy} \Big | \overrightarrow{x^{*}_{\epsilon Y} y_{l}} \Big |.
	\end{split}
\end{equation*}\end{linenomath}
Although a movement direction and a displacement magnitude can be this way determined for each boundary point  $\xi_{\epsilon Y}$, the actual relocation of  $\xi_{\epsilon Y}$ only occurs if sufficient but not complete ECM degradation will have occured in the peritumoural region $\epsilon Y \setminus \Omega(t_{0})$. To quantify the amount of ECM degradation, we use a \emph{transitional probability} that we define by
\begin{linenomath}\begin{equation*}
	q(x^{*}_{\epsilon Y}) := \dfrac{\int\limits_{\epsilon Y(x^{*}_{\epsilon Y}) \setminus \Omega(t_{0})} m(y, \tau) dy}{\int\limits_{\epsilon Y(x^{*}_{\epsilon Y})} m(y, \tau) dy},
\end{equation*}\end{linenomath}
Then, the movement of a boundary point is exercised only when adequate but not complete degradation of the peritunoural ECM occurs, which is characterized by the situation when this transitional probability $q(x^{*}_{\epsilon Y})$ exceeds a certain tissue threshold $\omega(\cdot, \cdot)$ (as defined in \cite{Dumitru_et_al_2013}), namely
\begin{linenomath}\begin{equation*}
	\omega(\beta,\! \epsilon Y)\! :=\!
	\begin{cases}
		\sin \! \bigg[ \dfrac{\pi}{2} \bigg( 1 \!-\! \dfrac{v(x^{*}_{\epsilon Y}, t_{0} + \Delta t)}{\beta \cdot \sup\limits_{\xi \in \partial \Omega(t_{0})} \!\! v(\xi, t_{0} \!+\! \Delta t)} \bigg) \bigg] & \!\! \text{if } \dfrac{v(x^{*}_{\epsilon Y}, t_{0} \!+\! \Delta t)}{\sup\limits_{\xi \in \partial \Omega(t_{0})} \!\! v(\xi, t_{0} \!+\! \Delta t)} \leq \beta, \\[30pt]
		\sin \! \bigg[ \dfrac{\pi}{2(1 \!-\! \beta)} \bigg( \dfrac{v(x^{*}_{\epsilon Y}, t_{0} + \Delta t)}{\sup\limits_{\xi \in \partial \Omega(t_{0})} \!\! v(\xi, t_{0} \!+\! \Delta t)} - \beta \bigg) \bigg] & \!\! \text{otherwise},
	\end{cases}
\end{equation*}\end{linenomath}
where $\beta \in (0, 1)$ controls the optimal level of ECM for cancer invasion and $v(x,t) := l(x,t) + F(x,t)$.

\end{paracol}
\reftitle{References}


\externalbibliography{yes}
\bibliography{ThesisReferences}

\end{document}